\documentclass[sn-mathphys-num]{sn-jnl}

\usepackage{graphicx}%
\usepackage{multirow}%
\usepackage{amsmath,amssymb,amsfonts}%
\usepackage{amsthm}%
\usepackage{mathrsfs}%
\usepackage{xcolor}%
\usepackage{textcomp}%
\usepackage{manyfoot}%
\usepackage{booktabs}%
\usepackage{algorithm}%
\usepackage{algorithmicx}%
\usepackage{algpseudocode}%
\usepackage{listings}%

\theoremstyle{thmstyleone}%
\newtheorem{theorem}{Theorem}

\newtheorem{proposition}[theorem]{Proposition}

\newtheorem{lemma}[theorem]{Lemma}

\theoremstyle{thmstyletwo}%
\newtheorem{example}{Example}%
\newtheorem{remark}{Remark}%

\theoremstyle{thmstylethree}%
\newtheorem{definition}{Definition}%

\def\Limsup{\mathop{{\rm Lim}\,{\rm sup}}}

\def\gph{\mbox{\rm gph}\,}
\def\epi{\mbox{\rm epi}\,}

\def\R{\mathbb{R}}
\def\N{\mathbb{N}}
\def\ox{\bar{x}}

\begin{document}
	
	\title[Second-order optimality conditions for $C^{1,1}$ sparsity constrained optimization problems]{Second-order optimality conditions for $C^{1,1}$ sparsity constrained optimization problems via disjunctive programming}

	\author[1,3]{\fnm{Luu} \sur{Thi Thu Huyen}}\email{luuthithuhuyen@hvu.edu.vn}
	\affil[1]{Faculty of Natural Sciences, Hung Vuong University, Phu Tho, Vietnam}
	
	\author[2]{\fnm{Liguo} \sur{Jiao}}\email{jiaolg356@nenu.edu.cn; hanchezi@163.com}
	\affil[2]{Academy for Advanced Interdisciplinary Studies, Northeast Normal University,
		Changchun 130024, Jilin Province, China}
		
	\author*[3]{\fnm{Nguyen} \sur{Van Tuyen}}\email{nguyenvantuyen83@hpu2.edu.vn}
	\affil*[3]{Department of Mathematics, Hanoi Pedagogical University 2, Xuan Hoa, Phu Tho, Vietnam}

	\abstract{In this paper, we investigate second-order necessary and sufficient optimality conditions for a class of $C^{1,1}$ sparsity constrained optimization problems via a disjunctive programming approach. More precisely, by utilizing the limiting/Mordukhovich second-order subdifferential of the associated Lagrangian function, we establish new second-order necessary and sufficient optimality conditions for local optimal solutions. These results are derived under mild assumptions and extend several existing results in the literature. In addition, we apply our theoretical developments to sparse multiobjective optimization problems and derive second-order optimality conditions for efficient solutions.  }


	
	\maketitle
	
\section{Introduction}
Sparsity constrained optimization has attracted considerable attention because of its applications in high-dimensional statistics, compressed sensing, signal processing, and machine learning. A typical sparsity constrained optimization problem can be formulated as
\begin{equation*} 
	\min f(x)\ \ \text{s.t.}\ \ x\in\Omega:=\{x\in \R^n\,:\|g(x)\|_0\leq s\}, 
\end{equation*}
where $f\colon\R^n\to\R$, $g\colon\R^n\to\R^m$,  are given functions,   $s\in (0, m)$ is a given positive integer,  and $\|\cdot\|_0$ denotes the number of nonzero components of a vector.  Such problems arise naturally in sparse regression, sparse approximation, and regularized inverse problems; see, e.g., \cite{CandesRombergTao2006,Donoho2006,Tibshirani1996,Hastie2009}. The sparsity constraint, however, leads to a highly nonconvex feasible set even when the functions involved are smooth, which makes the development of appropriate optimality conditions a challenging issue.
	
A useful feature of sparsity constraints is their disjunctive structure. In fact, the sparsity set
$$
	\mathcal{S}:=\{z\in\R^m\,:\, \|z\|_0\leq s\}
$$
	can be represented as a finite union of linear subspaces,
$$
	\mathcal{S}=\bigcup_{J\in\mathcal{J}}S_J,	
$$
where $\mathcal{J}:=\{J\subset\{1,\ldots,m\}\,:\,|J|=s\}.$ This representation makes it possible to treat sparsity constrained optimization within the framework of disjunctive programming. Movahedian et al. \cite{Movahedian-et-al} explicitly formulated sparsity constrained optimization problems as disjunctive programs and developed constraint qualifications and first- and second-order optimality conditions by exploiting the geometry of the sparsity set. Their work demonstrates the effectiveness of the disjunctive programming framework in dealing with the nonconvex geometry induced by sparsity constraints. More generally, optimality conditions and constraint qualifications for disjunctive programs have been extensively investigated; see, e.g., \cite{Flegel-et-al-07,Gfrerer2014,Benko-et-al-22,XuYe,Mehlitz}.
	
	Second-order optimality conditions provide substantially finer information about local solutions than first-order stationarity conditions and are fundamental in the analysis of stability and numerical algorithms. For smooth optimization problems, classical second-order conditions are expressed in terms of the Hessian of the Lagrangian. For nonsmooth problems, generalized second-order derivatives provide natural substitutes for the Hessian. In particular, second-order subdifferentials in the sense of variational analysis have been successfully employed to derive necessary and sufficient optimality conditions for nonsmooth optimization problems; see \cite{ChieuLeeYen2017,Huy-Tuyen-16,Huy-Tuyen-19,Khanh-etal-25,Mordukhovich_2006,Mordukhovich_2018,Mordukhovich_2024,Tuyen-Huy-Kim}.
	
	A relevant line of research concerns the use of second-order subdifferentials of Lagrangian functions. An et al. \cite{An-Xu-Yen} developed second-order necessary and sufficient optimality conditions for optimization problems with abstract constraints by using the Fr\'echet second-order subdifferential of the Lagrangian function. Their results show that second-order variational analysis can be applied even when the objective functions are only $C^1$-smooth. More recently, An and Tuyen \cite{An-Tuyen} investigated $C^{1,1}$ optimization problems with inequality and equality constraints by employing the limiting/Mordukhovich second-order subdifferential of the associated Lagrangian function. In particular, they established second-order necessary and sufficient optimality conditions under assumptions appropriate for $C^{1,1}$ data and considered applications to multiobjective optimization. Thus, the limiting second-order subdifferential provides a natural framework for second-order analysis beyond the twice continuously differentiable setting.
	
	The combination of second-order variational analysis with sparsity constraints has recently received further attention. Bucher and Schwartz \cite{BucherSchwartz2018} derived second-order optimality conditions and convergence results for regularization methods in cardinality-constrained optimization. Kan and Song \cite{Kan-Song} studied second-order conditions related to augmented Lagrange multipliers for sparse optimization. Chen et al. \cite{ChenJW-et_al_2024} investigated first- and second-order optimality conditions for nonsmooth sparsity multiobjective optimization problems. Most recently, Chen et al. \cite{ChenJW-et-at-2026} established second-order optimality conditions for sparse optimization by employing the Fr\'echet second-order subdifferential of the Lagrangian function. In particular, their necessary conditions are obtained under suitable Robinson-type constraint qualifications, while sufficient conditions are derived for sparse problems with equality constraints  \cite{ChenJW-et-at-2026}.
	
	Despite these developments, further investigation of second-order optimality conditions for sparsity constrained problems with $C^{1,1}$ data remains necessary. The use of the Fr\'echet second-order subdifferential, as in \cite{An-Xu-Yen,ChenJW-et-at-2026}, may be restrictive in situations where the corresponding Fr\'echet second-order subdifferential is empty; see, e.g., \cite[Example 3.1]{Huyen-Jiao-Tuyen}. In contrast, the limiting second-order subdifferential has more favorable existence properties for $C^{1,1}$ functions and is therefore particularly suitable for second-order analysis in this setting. At the same time, the results in \cite{An-Tuyen} concern $C^{1,1}$ optimization problems with classical inequality and equality constraints rather than the highly nonconvex disjunctive structure generated by sparsity constraints. These observations motivate us to develop a second-order theory that combines the $C^{1,1}$ framework with the disjunctive geometry of sparsity constraints.
\medskip	

The purpose of this paper is to establish second-order necessary and sufficient optimality conditions for a class of $C^{1,1}$ sparsity constrained optimization problems via a disjunctive programming approach. More precisely, we exploit the representation of the sparsity set as a finite union of linear subspaces and employ the limiting second-order subdifferential of the associated Lagrangian function. This enables us to formulate second-order optimality conditions that are applicable to $C^{1,1}$ objective and constraint functions without requiring twice continuous differentiability.
\medskip
	
Our main contributions can be summarized as follows: 
\begin{itemize}
\item First, we establish second-order necessary optimality conditions for local optimal solutions of the considered sparsity constrained problem in terms of the limiting second-order subdifferential of the Lagrangian function. The derivation takes advantage of the disjunctive structure of the sparsity set and the associated normal and tangent cone constructions. 
\item Second, we  derive  second-order sufficient optimality conditions under appropriate assumptions, including conditions guaranteeing isolated local solutions of order $2$. These results complement the disjunctive programming results of \cite{Movahedian-et-al} and the Fr\'echet second-order subdifferential approach of \cite{ChenJW-et-at-2026}, while extending the $C^{1,1}$ limiting second-order framework developed in \cite{An-Tuyen} to sparsity constrained problems. In particular, the use of the limiting second-order subdifferential allows us to treat cases in which the corresponding Fr\'echet second-order subdifferential may fail to provide a suitable second-order object.
\item Third, we apply the obtained results to sparse multiobjective optimization problems. By combining scalarization with the proposed second-order analysis, we derive second-order necessary and sufficient optimality conditions for efficient solutions. This extends the scope of the proposed approach from scalar sparsity constrained optimization to multiobjective optimization with sparsity constraints.
\end{itemize}
	 
It is worth mentioning the recent work of Huyen et al. \cite{Huyen-Jiao-Tuyen}, which also employs the limiting   second-order subdifferential to derive second-order optimality conditions for $C^{1,1}$ sparse optimization problems. Their approach formulates the feasible set as the intersection of a conventional constraint system and the sparsity set, and relies on restricted constraint qualifications together with explicit tangent and normal cone formulas for the sparsity set. In contrast, the present paper adopts a disjunctive programming viewpoint by representing the sparsity constraint as a finite union of linear subspaces and exploiting the associated disjunctive structure directly. This leads to a different treatment of the active subspaces and the corresponding first- and second-order geometric objects. In particular, our analysis allows the second-order optimality conditions to be formulated in terms of the geometry of the active disjunctive branches and the limiting second-order subdifferential of the associated Lagrangian. Thus, although both works employ the limiting second-order subdifferential in the $C^{1,1}$ setting, the underlying variational frameworks and the resulting optimality conditions are different.
\medskip
	
The remainder of the paper is organized as follows. Section \ref{Section2} recalls some concepts and results from variational analysis  that will be used throughout the paper. Section \ref{Section 3} is devoted to establishing second-order necessary and sufficient optimality conditions for the considered sparsity constrained problem. Section \ref{Section 4}    applies the developed theory to sparse multiobjective optimization problems. Concluding remarks are given in the final section.

\section{Preliminaries} \label{Section2}
Our notation and terminology are standard in variational analysis; see, e.g., the books \cite{Aubin1990,Borwein2005,Mordukhovich_2006,Mordukhovich_2018,Mordukhovich_2024,Rockafellar1998}. Recall that for any number $n\in\mathbb{N}:=\{1,2,\ldots\}$ we denote $x:=(x_1,\ldots,x_n)$ and equip the space $\mathbb{R}^n$ with the usual scalar product $\langle\cdot,\cdot\rangle$ and the Euclidean norm $\|\cdot\|$.  For a nonempty set $X \subset \mathbb{R}^n,$ the closure,  convex hull and conic hull of $X$ are denoted, respectively, by $\mathrm{cl}\,X$,  $ \mathrm{co}\,X$, and $\mathrm{cone}\,X.$ If $X$ is a cone, then  the polar cone of $X$, denoted by $X^\circ$, is defined by
\begin{equation*}
	X^\circ:=\{u\in\R^n\,:\, \langle u, x\rangle\leq 0 \ \ \forall x\in X\}.
\end{equation*} 

Let  $e_i$, $i=1, \ldots, n$,  be the {\em  $i$-th unit vector} of $\R^n$. For a given subset $J\subset\{1, \ldots, n\}$, $\text{span}\,\{e_i\,:\, i\in J\}$  denotes the subspace of $\R^n$ spanned by the vectors $\{e_i\,:\, i\in J\}$.  

For a set-valued mapping $F \colon \mathbb{R}^n \rightrightarrows \mathbb{R}^m$, the {\em Painlev\'e--Kuratowski outer limit} of $F$ at $\bar x\in\R^n$ is defined by
\begin{equation*}
	\Limsup_{x \to {\bar x}} F(x) := \{y \in \mathbb{R}^m \,:\, \exists x_k \to {\bar x}, \exists y_k \in F(x_k), y_k \to y\}.
\end{equation*}

\begin{definition}\rm 
	Let $X\subset \R^n$ and $\bar x\in X$. The set
	\begin{equation*}
		T_X(\bar x):=\{d\in\R^n\,:\, \exists t_k\downarrow 0, \exists d_k\to d, \bar x+t_kd_k\in X \ \ \forall k\in\N\}
	\end{equation*}
	is called the {\em contingent cone} (or the {\em Bouligand--Severi
		tangent cone}) to $X$ at $\bar x$. 
\end{definition}

It is well-known that $T_X(\bar x)$ is a nonempty closed cone and when $X$ is a convex  set, then
\begin{equation*}
	T_X(\bar x)=\mathrm{cl}\,\mathrm{cone}(X-\bar x).
\end{equation*}

\begin{definition}[{see~\cite[pp. 5--6]{Mordukhovich_2018}}]
	\rm	Let $X$ be a nonempty subset of $\mathbb{R}^n$ and $\bar x \in X$. The \textit{Fr\'echet/regular normal cone}  to $X$ at $\bar x$ is defined by
	\begin{align*}
		\widehat N_X(\bar x):=\Big\{ v\in \mathbb{R}^n \,:\, \limsup\limits_{x \xrightarrow{X}\bar x} \dfrac{\langle v, x-\bar x \rangle}{\|x-\bar x\|} \leq 0 \Big\},
	\end{align*}
	where $x \xrightarrow{X} \bar x$ means that $x \rightarrow \bar x$ and $ x\in X$. The \textit{limiting/Mordukhovich normal cone} to $X$ at $\bar x$ is given by
	\begin{align*}
		N_X(\bar x):=\Limsup\limits_{ x \xrightarrow{X} \bar x} \widehat{N}_X(x).
	\end{align*}
	We put $\widehat N_X(\bar x)=N_X(\bar x) :=\emptyset$ if $\bar x \not\in X$.
\end{definition}	

By definition, it is clear that $\widehat N_X(\bar x) \subset N_X(\bar x)$  for all $\bar x \in X.$ Furthermore, it is well-known that $\widehat N_X(\bar x)=[T_X(\bar x)]^\circ$   and it is always convex. When $X$ is a convex set, then the Fr\'echet normal cone, the limiting normal cone, and the normal cone in the sense of convex analysis coincide, i.e., 
\begin{align*}
	\widehat N_X(\bar x) = N_X(\bar x):=\{ v\in \mathbb{R}^n : \langle v,x- \bar x \rangle \le 0,  \ \forall x\in X\}.
\end{align*}

Let  $F: \mathbb{R}^n\rightrightarrows{\mathbb{R}^m}$ be a  multifunction with the \textit{domain} $${\rm{dom}}\, F:=\{ x \in \mathbb{R}^n : F(x)\not=\emptyset\}$$ and the \textit{graph} $${\rm{gph}}\, F:=\{ (x,y) \in \mathbb{R}^n \times \mathbb{R}^m : y \in F(x)\}.$$  

\begin{definition} [{see~\cite[Definition 1.11]{Mordukhovich_2018}}]\rm  Let $  F: \mathbb{R}^n\rightrightarrows{\mathbb{R}^m}$ be a multifunction and  $(\bar x, \bar y) \in  {\rm{gph}}\, F$.	The  \textit{limiting/Mordukhovich coderivative}  of  $F$ at $(\bar x, \bar y)$ is a multifunction $D^* F(\bar x, \bar y): \mathbb{R}^m \rightrightarrows{\mathbb{R}^n}$ with the values
	\begin{align*}
		D^* F(\bar x,\bar y)(v):=&\left\{u\in \mathbb{R}^n : (u, -v) \in  N_{\gph F}(\bar x, \bar y)\right\}, \quad v \in \mathbb{R}^m.
	\end{align*}
	If $(\bar x, \bar y) \notin {\rm{gph}}\, F$, one puts $ D^* F(\bar x, \bar y)(v):=\emptyset$ for any $v\in \mathbb{R}^m$. The symbol $D^* F(\bar x)$ is used when $F$ is single-valued at $\bar x$ and $\bar y=F(\bar x)$. 
\end{definition}

Consider a function $\varphi: \mathbb{R}^n\rightarrow \overline{\mathbb{R}}=\mathbb{R}\cup \{+\infty\}$ with the \textit{effective domain} 
$$\mbox{dom}\, \varphi:=\{x \in \mathbb{R}^n : \varphi(x) < +\infty\}$$
and the \textit{epigraph} 
$$ {\rm{epi}}\, \varphi:=\{ (x, \alpha) \in \mathbb{R}^n \times \mathbb{R} : \varphi (x)\leq \alpha\}.$$

\begin{definition}[{see~\cite[Definition 1.77]{Mordukhovich_2006}}] \rm  
	Let $\bar {x} \in \mbox{dom}\, \varphi$. The  \textit{limiting/Mordukhovich subdifferential} of $\varphi$ at $\bar {x}$ is defined by
	\begin{align*}
		\partial \varphi(\bar x):=\left\{v \in \mathbb{R}^n : (v, -1) \in  N_{\epi \varphi}(\bar {x}, \varphi(\bar{ x}))\right\}.
	\end{align*}
	If $\bar {x} \notin \mbox{dom}\, \varphi$, then we put  $\partial \varphi (\bar x):=\emptyset$.
\end{definition}

One can use the notion of coderivative to construct the second-order generalized differential theory of extended-real-valued functions.

\begin{definition}  [{see \cite[Definition 3.17]{Mordukhovich_2018}}] {\rm Let $\varphi:\mathbb{R}^n\rightarrow \overline{\mathbb{R}}$ 	be a function and  $\bar{x}\in\mathrm{dom}\,\varphi.$ For any $\bar y\in  \partial \varphi (\bar x)$, the map $\partial^2 \varphi(\bar x,\bar y):\mathbb{R}^n\rightrightarrows \mathbb{R}^n$ with the values
		\begin{align*}
			\partial^2 \varphi(\bar x,\bar y)(v):=(D^* \partial \varphi)(\bar x,\bar y)(v)=\{u : (u,-v)\in N_{\gph \partial \varphi}(\bar x, \bar y)\}
		\end{align*} is said to be the \textit{limiting/Mordukhovich
			second-order subdifferential} of $\varphi$ at $\bar x$ relative to	$\bar y.$}
\end{definition}

If $\partial \varphi (\bar x)$ is a singleton, the symbol $\bar y$ in the notation $\partial^2 \varphi (\bar x,\bar y)(v)$ will be omitted.  When $\varphi$ is a $C^2$-smooth function around $\bar x$, i.e., $\varphi$ is twice continuously differentiable in a neighborhood of $\ox$, then  
\begin{align*}
	\partial^2 \varphi(\bar x)(v)=\{\nabla^2 \varphi (\ox)^* v\}= \{\nabla^2 \varphi(\bar x)^T v\} \ \  v\in \mathbb{R}^n
\end{align*}
with $\nabla^2 \varphi (\ox)$ being the Hessian matrix of $\varphi$ at $\ox$ and the apex $^T$ denoting the matrix transposition; see, e.g., \cite[p. 124]{Mordukhovich_2018}.

\medskip
Let $U\subset\R^n$ be a nonempty, open, and convex set. We denote by  $C^{1,1}(U)$ the class of all real-valued functions $\varphi$, which are Fr\'echet differentiable on $U$, and whose gradient mapping $\nabla \varphi (\cdot)$ is locally Lipschitz on $U$. By \cite[Theorem 1.90]{Mordukhovich_2006}, if $\varphi\in C^{1,1}(U)$ and $\bar x\in U$, then one has
\begin{align*} 
	\partial^2 \varphi(\bar x) (v) = \partial^2 \varphi (\bar x, \nabla \varphi (\bar x))(v)=\partial \langle v, \nabla \varphi \rangle (\bar x)\ \ \forall v\in \mathbb{R}^n.
\end{align*}

The following properties can be obtained directly from the definition.
\begin{proposition}\label{property}
	Let $\varphi$ be of class $C^{1,1}(U)$ and $\bar x\in U$. The following assertions hold:
	\begin{enumerate}
		\item[\rm (i)] $\partial^2 \varphi(\bar x)(\lambda v) = \lambda \partial^2 \varphi(\bar x)(v)$ for all $\lambda \ge 0$ and $v\in \mathbb{R}^n.$
		\item[\rm (ii)] For any $v\in \mathbb{R}^n$ the mapping $x \mapsto \partial^2 \varphi(x) (v)$ is locally bounded. Moreover, if $x_k \to \bar x$, $x_k^* \to x^*$, $x_k^* \in \partial^2 \varphi(x_k)(v)$ for all $k\in \mathbb{N}$, then $x^*\in \partial^2 \varphi(\bar x)(v).$
	\end{enumerate} 
\end{proposition}

The Taylor formula for $C^{1,1}$ functions in terms of the limiting second-order subdifferential plays a crucial role in our investigation.
\begin{theorem} [{see~\cite[Theorem 3.1]{Feng_Li_2020}}]\label{Taylor_formula}
	Let $\varphi$ be of class $C^{1,1}(U)$  and $a, b\in U$. Then, there exist $z\in \partial^2 \varphi(x) (b-a)$ and $z'\in \partial^2 \varphi(x') (b-a)$ for some $x$ and $x' \in [a,b]$  such~that
	$$ \dfrac{1}{2} \langle z', b-a \rangle \le \varphi(b)-\varphi(a)-\langle\nabla \varphi(a), b-a \rangle \le \dfrac{1}{2}\langle z, b-a \rangle.$$
\end{theorem}

\section{Second-order optimality conditions}\label{Section 3}

Consider the following sparsity constrained optimization problem
\begin{equation}\label{Problem}
	\min f(x)\ \ \text{s.t.}\ \ x\in\Omega:=\{x\in \R^n\,:\|g(x)\|_0\leq s\},\tag{SP}	
\end{equation}
where $f\colon\R^n\to\R$, $g\colon\R^n\to\R^m$,  are $C^{1,1}$ functions,   $s\in (0, m)$ is a given positive integer,  and $\|\cdot\|_0$ denotes the number of nonzero components of a vector.

For a subset $J$ of $\{1, \ldots, m\}$, we denote the cardinality and the complement  of    $J$, respectively, by $|J|$ and $J^c$. Put $$\mathcal{S}:=\{z\in\R^m\,:\, \|z\|_0\leq s\},$$  
$$\mathcal{J}:=\{J\subset \{1, \ldots, m\}\,:\, |J|=s\},$$ 
and  $S_J:=\mathrm{span}\{e_j\,:\, j\in J\}$. Then, one has
\begin{equation*}
	\mathcal{S}=\bigcup_{J\in \mathcal{J}} S_J,
\end{equation*} 
and 
\begin{equation*}
	\Omega=\{x\in\R^n\,:\, g(x)\in\mathcal{S}\}.
\end{equation*}
For $\bar x\in \Omega$, we denote 
$$\mathcal{J}(\bar x):=\{J\in\mathcal{J}\,:\, g(\bar x)\in S_J\}.$$  

The support index set of a point $z\in\R^m$, denoted by $\Gamma(z)$, is defined by
\begin{equation*}
	\Gamma(z):=\{i\in \{1, \ldots, m\}\,:\, z_i\neq 0\}.
\end{equation*} 

The following result provides formulae for computing the contingent cone, the Fr\'echet normal cone, and the Mordukhovich normal cone to $\mathcal{S}$. 
\begin{lemma}[{see \cite[Lemma 3.1]{Movahedian-et-al}}]\label{Lemma3} For any $\bar x\in\Omega,$ we have
	\begin{align*}
		T_\mathcal{S}(g(\bar x))=\bigcup_{J\in \mathcal{J}(\bar x)} S_J,
	\end{align*}	

	\begin{equation*}
		\widehat{N}_\mathcal{S}(g(\bar x))=
		\begin{cases}
			\mathrm{span}\{e_i\,:\, i\notin \Gamma(g(\bar x))\}, \ \ &\text{if}\ \ |\Gamma(g(\bar x))|=s,
			\\
			\{0_{\R^m}\}, \ \ &\text{if}\ \ |\Gamma(g(\bar x))|<s,
		\end{cases}
	\end{equation*}
and 
\begin{equation*}
N_\mathcal{S}(g(\bar x))=\bigcup_{J\in \mathcal{J}(\bar x)} S_{J^c}.
\end{equation*}
\end{lemma}

The Lagrangian function $\mathcal{L}\colon\R^n\times \R^m\to \R$ associated with the constrained problem \eqref{Problem} is given by
$$\mathcal{L}(x, \lambda):=f(x)+\lambda^T g(x)$$
for each $(x, \lambda)\in \R^n\times \R^m$. The set of Lagrange multipliers at $\bar x\in\Omega$ is defined by
\begin{equation*}
	\Lambda(\bar x):=\{\lambda\in\R^m\,:\, \exists J\in \mathcal{J}(\bar x) \ \ \text{such that}\ \  \lambda \in S_{J^c}, \nabla_x\mathcal{L}(\bar x, \lambda)=0\}.
\end{equation*}

For  $\bar x\in\Omega$, the {\em linearized cone} $L(\bar x)$ to $\Omega$ at $\bar x$ is defined by
\begin{equation*}
	L(\bar x):=\{d\in\R^n\,:\, \nabla g(\bar x)d\in T_\mathcal{S}(g(\bar x))\}.
\end{equation*} 

The following lemma provides an explicit characterization of $L(\bar x)$ and its negative polar cone.
\begin{lemma}[{see \cite[Lemma 3.3]{Movahedian-et-al}}]\label{Lemma4}  For  $\bar x\in\Omega,$ we have
	\begin{align*}
		&L(\bar x)=\bigcup_{J\in\mathcal{J}(\bar x)}\{d\in\R^n\,:\,\langle \nabla g_i(\bar x), d\rangle =0 \ \ \forall i\in J^c\} \ \ \text{and}
		\\
		&[L(\bar x)]^\circ=\bigcap_{J\in \mathcal{J}(\bar x)}\mathrm{span}\{\nabla g_i(\bar x)\,:\, i\in J^c\}.
	\end{align*}	
\end{lemma}

Next, we recall several notions of stationarity for \eqref{Problem}  based on the gradient vectors of the objective and constraint functions.

\begin{definition}[{see \cite[Definition 3.1]{Movahedian-et-al}}]\rm Let $\bar x\in\Omega$. We say that
	\begin{enumerate}
		\item[\rm(i)] $\bar x$ is an {\em  $M$-stationary point} of \eqref{Problem} if
		\begin{equation*}
			0\in\nabla f(\bar x)+\bigcup_{J\in \mathcal{J}(\bar x)}\mathrm{span}\{\nabla g_i(\bar x)\,:\, i\in J^c\}.
		\end{equation*}
		\item[\rm(ii)] $\bar x$ is an {\em  $L$-stationary point} of \eqref{Problem} if
		\begin{equation*}
			0\in\nabla f(\bar x)+\bigcap_{J\in \mathcal{J}(\bar x)}\mathrm{span}\{\nabla g_i(\bar x)\,:\, i\in J^c\}.
		\end{equation*}
		\item[\rm(iii)] $\bar x$ is an {\em  $S$-stationary point} of \eqref{Problem} if
		\begin{equation*}
			\begin{cases}
				0\in\nabla f(\bar x)+\mathrm{span}\{\nabla g_i(\bar x)\,:\, i\notin \Gamma(g(\bar x))\},  \ \ &\text{if}\ \ |\Gamma(g(\bar x))|=s,
				\\
				0=\nabla f(\bar x), \ \ &\text{otherwise}.
			\end{cases}
		\end{equation*}
	\end{enumerate}
\end{definition}

By definition, it is clear that 
\\
{\centerline{$S$-stationarity $\Rightarrow$ $L$-stationarity $\Rightarrow$ $M$-stationarity}}
\\
and when $|\Gamma(g(\bar x))|=s$, then all the above stationary notions coincide. We note that the terms $M$-stationary and $S$-stationary are derived from the formulae for the Mordukhovich and Fr\'echet normal cones to $\mathcal{S}$ given in Lemma \ref{Lemma3}, respectively, whereas the term $L$-stationary is derived from the formula for the polar cone of $L(\bar x)$ given in Lemma \ref{Lemma4}; see \cite{Flegel-et-al-07,Movahedian-et-al}.

\medskip

Next, we recall some notions of optimal solutions to the problem \eqref{Problem}.

\begin{definition}\rm    Let $\bar x\in\Omega$. We say that:
	\begin{enumerate}
		\item[\rm(i)] $\bar x$ is a {\em  local optimal solution} to the problem \eqref{Problem} if there exists a neighborhood $U$ of $\bar x$ such that $f(x)\geq f(\bar x)$ for all $x\in\Omega\cap U$. 
		
		\item[\rm(ii)]  $\bar x$ is a {\em  strict local  optimal solution} to the problem \eqref{Problem} if there exists a neighborhood $U$ of $\bar x$ such that $f(x)> f(\bar x)$ for all $x\in[\Omega\cap U]\setminus\{\bar x\}$.
		
		\item[\rm(iii)] (see \cite[Definition 2.1]{Auslender_1984}) $\bar x$ is an \textit{isolated local solution of order 2}  to the problem \eqref{Problem} if there exist a  constant $\rho>0$ and a neighborhood $U$ of $\bar x$ such that
		\begin{equation*} 
			f (x)\geq  f (\bar x)  +\rho \|x-\bar x\|^2\ \ \forall x\in \Omega \cap U.
		\end{equation*}
	\end{enumerate}	
\end{definition}

Clearly, every isolated local solution of order $2$ to the problem \eqref{Problem} is a strict local optimal solution, and every strict local optimal solution to the problem \eqref{Problem} is a local optimal solution. The converses of these implications do not hold in general. For example, if $f(x)$ is constant on $\mathbb{R}$, then every $\bar x\in\mathbb{R}$ is a local optimal solution, but none is a strict local optimal solution. Furthermore, for $f(x)=x^4$ for all $x\in\mathbb{R}$, then $\bar x=0$ is a strict local optimal solution of $f$ on $\mathbb{R}$, but it is not an isolated local solution of order $2$.

\medskip     

We recall the following constraint qualifications for \eqref{Problem}, which were introduced in \cite{Movahedian-et-al}; see also \cite{Flegel-et-al-07}.  

\begin{definition} [{see \cite[Definition 3.3]{Movahedian-et-al}}]\rm Let $\bar x\in\Omega$. We say that $\bar x$ satisfies  
	\begin{enumerate}
		\item [\rm(i)] the {\em Abadie constraint qualification} (ACQ) if $T_\Omega(\bar x)= L(\bar x)$.
		\item [\rm(ii)] the {\em Guignard constraint qualification} (GCQ) if $\widehat{N}_\Omega(\bar x)=[L(\bar x)]^\circ$.
		
		\item [\rm(iii)] the {\em  Mangasarian Fromovitz constraint qualification} (MFCQ) if for each $J\in\mathcal{J}(\bar x)$, the system $\{\nabla g_i(\bar x)\,:\, i\in J^c\}$ is  linear  independent.
	\end{enumerate}
\end{definition}

\begin{remark}\rm It follows from \cite{Flegel-et-al-07} that the following implications hold:

		\centering{ MFCQ $\Rightarrow$ ACQ $\Rightarrow$ GCQ.}
\end{remark}

\medskip

The next theorem provides first-order optimality conditions for a local optimal solution of \eqref{Problem} under the weakest constraint qualification, namely, the GCQ.

\begin{lemma}[{see \cite[Theorem 3.1]{Movahedian-et-al}}]\label{Lemma31} Let $\bar x\in\Omega$ be a local optimal solution to the problem \eqref{Problem}. Then the following assertions hold:
\begin{enumerate} 
	\item [\rm(i)] If $\bar x$ satisfies the GCQ, then $\bar x$ is an $L$-stationary point of \eqref{Problem}.  
	
	\item  [\rm(ii)] If $\bar x$ satisfies the GCQ and the following intersection property
	\begin{equation}\label{New-equa-10}
	\|g(\bar x)\|_0<s \Rightarrow  \bigcap_{J\in \mathcal{J}(\bar x)}\mathrm{span}\{\nabla g_i(\bar x)\,:\, i\in J^c\}=\{0\},
	\end{equation} 
then $\bar x$ is an $S$-stationary point of \eqref{Problem}. 
\end{enumerate}
 
\end{lemma}

We next introduce several types of second-order necessary optimality conditions for \eqref{Problem} in terms of limiting second-order subdifferentials.  

\begin{definition}\rm Let $\bar x$ be an $M$-, an $L$- or an $S$-stationary point of \eqref{Problem}. We say that:
	\begin{enumerate}
		\item [\rm(i)] the {\em basic second-order necessary optimality condition} (BSOC) (resp., the {\em strict basic second-order necessary optimality condition} (SBSOC)) holds at $\bar x$, if for each $d\in L(\bar x)$ (resp., $d\in L(\bar x)\setminus\{0\}$) there exist  $\lambda\in\Lambda(\bar x)$ and $z\in \partial^2 \mathcal{L}(\bar x, \lambda)(d)$ such that
		\begin{equation*}
			\langle z, d\rangle \geq 0 \ \ (\text{resp.,}\ \ \langle z, d\rangle > 0).	
		\end{equation*}   
		
		\item [\rm(ii)] the {\em strong second-order necessary optimality condition} (SSOC) (resp., the {\em strict strong second-order necessary optimality condition} (SSSOC)) holds at $\bar x$, if there exists  $\lambda\in\Lambda(\bar x)$ such that, for each $d\in L(\bar x)$ (resp., $d\in L(\bar x)\setminus\{0\}$), there exists $z\in \partial^2 \mathcal{L}(\bar x, \lambda)(d)$ satisfying
		\begin{equation*}
			\langle z, d\rangle \geq 0 \ \ (\text{resp.,}\ \ \langle z, d\rangle > 0).	
		\end{equation*}   
	\end{enumerate}	 
\end{definition}

\begin{remark}\rm By definition, the following implications hold:
\begin{enumerate}
\item[\rm(i)] SBSOC $\Rightarrow$ BSOC;

\item[\rm(ii)] SSSOC  $\Rightarrow$ SSOC $\Rightarrow$ BSOC.
\end{enumerate}	
Furthermore, we observe that a drawback of the BSOC (resp., SBSOC), compared with the SSOC (resp., SSSOC), is that the Lagrange multiplier $\lambda$ need not be the same for all critical vectors. Hence, from a computational point of view, the SSOC (resp., SSSOC) is more meaningful than the BSOC (resp., SBSOC).
\end{remark}

The following theorem provides second-order necessary optimality conditions for \eqref{Problem} in terms of the limiting second-order subdifferential of the associated Lagrangian function.

\begin{theorem}[Necessary optimality conditions]\label{Necessary-condition} Let $\bar x\in\Omega$ be a feasible solution to the problem \eqref{Problem} which satisfies the ACQ. 
	\begin{enumerate}
		\item[\rm(i)] If $\bar x$ is a local optimal solution to the problem \eqref{Problem}$,$ then the BSOC holds at $\bar x$. In addition$,$ if the {\em  locally complementary slackness condition} holds at $\bar x,$ i.e.$,$ there exists a neighborhood $U$ of $\bar x$ and $\bar J\in\mathcal{J}(\bar x)$ such that $g(x)\in S_{\bar J}$ for all $x\in \Omega\cap U,$ then the SSOC holds at $\bar x$.  
		
		\item[\rm(ii)] If $\bar x$ is an isolated local solution of order $2$ to the problem \eqref{Problem}$,$ then the SBSOC holds at $\bar x.$ In addition$,$ if the locally complementary slackness condition holds at $\bar x,$ then the SSSOC holds at $\bar x$.    
	\end{enumerate}	
\end{theorem}
\begin{proof} (i) We first prove that if $\bar x$ is a local optimal solution to the problem \eqref{Problem}, then the BSOC holds at $\bar x$. Indeed, since $\bar x$ is a local optimal solution to the problem \eqref{Problem} and the ACQ holds at $\bar x$, we see that $\bar x$ is an $L$-stationary point and there exists a neighborhood $U$ of $\bar x$ such that
	\begin{equation}\label{New-equa-2}
		f(x)\geq f(\bar x) \ \ \forall x\in \Omega\cap U.
	\end{equation}	
	
	Fix any $d\in L(\bar x)$, then by the ACQ, one has $d\in T_\Omega(\bar x)$. By definition, there exist sequences $t_k\downarrow 0$ and $d_k\to d$ such that $x_k:=\bar x+t_kd_k\in\Omega\cap U$ for all $k\in\N$. Clearly, $x_k\to \bar x$ as $k\to\infty$, and by the continuity of $g$, one has $g(x_k)\to g(\bar x)$. This, together with the fact that $g(x_k)\in\mathcal{S}$ for all $k\in\N$, implies that there exists $k_0\in\N$ such that $\Gamma(g(\bar x))\subset \Gamma (g(x_k))$ for all $k\geq k_0$. For each $k\geq k_0$, since $\Gamma(g(\bar x))\subset \Gamma (g(x_k))$ and $|\Gamma (g(x_k))|\leq s$, there exists $J_k\in\mathcal{J}(\bar x)$ such that $\Gamma (g(x_k))\subset J_k$. Hence $J_k\in\mathcal{J}(\bar x)$  for all $k\geq k_0$. Since $\mathcal{J}(\bar x)$ is finite, by passing to a subsequence if necessary, we may assume that there exists $\bar J\in \mathcal{J}(\bar x)$ such that $\Gamma (g(x_k))\subset \bar J$ for all $k\in\N$. Clearly, $\Gamma(g(x_k)-g(\bar x))\subset \Gamma(g(x_k))\cup \Gamma(g(\bar x))$ and hence we may assume that
	\begin{equation}\label{New-equa-11} 
		\Gamma(g(\bar x))\cup \Gamma(g(x_k)-g(\bar x))\subset \Gamma(g(x_k)) \subset \bar J\in\mathcal{J}(\bar x) \  \ \forall k\in\N.	
	\end{equation} 
	This, together with the fact that $\bar x$ is an $L$-stationary point, i.e,
	\begin{equation*}
		0\in\nabla f(\bar x)+\bigcap_{J\in\mathcal{J}(\bar x)}\mathrm{span}\{\nabla g_i(\bar x)\,:\, i\in J^c\},
	\end{equation*} 
	implies that there exists $\bar \lambda\in S_{\bar J^c}$ such that $\nabla f(\bar x)+\sum_{i\in\bar J^c}\bar\lambda_i\nabla g_i(\bar x)=0.$ Hence $\bar\lambda_i=0$ for all $i\in\bar J$, $\nabla_x\mathcal{L}(\bar x, \bar \lambda)=0$, and $\langle \bar\lambda, g(x_k)-g(\bar x)\rangle=0$ for all $k\in\N$. This and \eqref{New-equa-2} imply that
	\begin{equation*}
		\mathcal{L} (x_k, \bar \lambda)- \mathcal{L} (\bar x, \bar \lambda)= f(x_k)-f(\bar x) \ge 0 \ \ \forall k\in\N.
	\end{equation*}
	Hence for each $k\in\N$ one has
	\begin{align}\label{Taylor_expres_2n} \nonumber
		0&\le 	\mathcal{L} (x_k, \bar \lambda)- \mathcal{L} (\bar x, \bar \lambda)\\ \nonumber
		&= [ \mathcal{L} (x_k, \bar \lambda) -\mathcal{L} (\bar x + t_k d, \bar \lambda)] +[\mathcal{L} (\bar x + t_k d, \bar \lambda)- \mathcal{L} (\bar x, \bar \lambda) ]\\ 
		& = [ \mathcal{L} (x_k, \bar \lambda) -\mathcal{L} (\bar x + t_k d, \bar \lambda)] +[\mathcal{L} (\bar x + t_k d, \bar \lambda)- \mathcal{L} (\bar x, \bar \lambda) -  \langle \nabla_x \mathcal{L}(\bar x, \bar \lambda) , t_k d \rangle].
	\end{align} 
	In view of the assumptions, it follows that $\mathcal{L} (\cdot, \bar \lambda) \in C^{1,1} (\mathbb{R}^n)$. Let $L>0$ be a Lipschitz constant of $\nabla_x \mathcal{L} (\cdot, \bar \lambda)$. By the classical mean value theorem (see, e.g., \cite[Theorem~5.10]{Rudin_PMA1976}), there exists $u_k\in (\bar x+t_k d, x_k)$ such that
	\begin{equation*}
		\mathcal{L}_{1k}:=	\mathcal{L}(x_k,\bar\lambda)  -\mathcal{L}(\bar x+t_k d, \bar\lambda) = \langle \nabla_x \mathcal{L}(u_k,\bar\lambda), t_k(d_k-d) \rangle.
	\end{equation*}
	Hence
	\begin{align*}\nonumber
		|\mathcal{L}_{1k}|&= |\langle \nabla_x \mathcal{L}(u_k,\bar\lambda), t_k(d_k-d) \rangle| \\ \nonumber
		&= |\langle \nabla_x \mathcal{L}(u_k,\bar\lambda) -   \nabla_x \mathcal{L}(\bar x,\bar\lambda), t_k(d_k -d) \rangle| \\ \nonumber
		& \le t_k \| \nabla_x \mathcal{L}(u_k,\bar\lambda)-   \nabla_x \mathcal{L}(\bar x,\bar\lambda)\|  \cdot \|d_k-d\|
		\\
		&\le L t_k \|u_k-\bar x\| \cdot \|d_k-d\|.
	\end{align*}
	This, together with the fact that 
	\begin{align*}
		\|u_k-\bar x\| &\le \|u_k-(\bar x+t_kd)\|+t_k\|d\|
		\\
		&\leq t_k(\|d_k-d\|+\|d\|),  
	\end{align*}
	implies that
	\begin{equation*} 
		|\mathcal{L}_{1k}|\le L t_k^2(\|d_k-d\|+\|d\|)\|d_k-d\|.
	\end{equation*}
	Hence
	\begin{equation}\label{New-equa-7} 
		\dfrac{|\mathcal{L}_{1k}|}{t_k^2} \le L (\|d_k-d\|+\|d\|)\|d_k-d\| \to 0 \ \ \text{as}\ \ k\to\infty.
	\end{equation}
	
	On the other hand, for each $k\in\N$, by applying the Taylor formula (Theorem~\ref{Taylor_formula}) for the function $\mathcal{L}(\cdot, \bar\lambda)$ on $[\bar x, \bar x+t_k d]$,  we can find  $w_k\in \partial^2 \mathcal{L} (\xi_k, \bar \lambda) (t_k d)$ for some $\xi_k\in [\bar x, \bar x+t_k d]$ such that 
	\begin{align*} 
		\mathcal{L}_{2k}:=\mathcal{L} (\bar x + t_k d, \bar \lambda)- \mathcal{L} (\bar x, \bar \lambda) -  \langle \nabla_x \mathcal{L}(\bar x, \bar \lambda) , t_k d \rangle \le  \dfrac{1}{2} \langle w_k, t_k d\rangle.
	\end{align*}
	By Proposition~\ref{property}(i), $w_k=t_kz_k$ for some $z_k\in \partial^2 \mathcal{L} (\xi_k, \bar \lambda)(d).$  Thus
	\begin{align}\label{Taylor_expres_4}
		\mathcal{L}_{2k} \le  \dfrac{1}{2}t_k^2 \langle z_k, d\rangle.
	\end{align}
	Combining \eqref{Taylor_expres_2n}  together with \eqref{Taylor_expres_4} we obtain
	\begin{align*}
		0\le \mathcal{L}_{1k} +\mathcal{L}_{2k} \le  \dfrac{1}{2}t_k^2 \langle z_k, d\rangle+\mathcal{L}_{1k},
	\end{align*}
	or, equivalently,
	\begin{equation}\label{New-equa-1}
		0\le  \dfrac{1}{2} \langle z_k, d\rangle+\dfrac{\mathcal{L}_{1k}}{t_k^2}.
	\end{equation}
	Since $\mathcal{L} (\cdot, \bar \lambda) \in C^{1,1} (\mathbb{R}^n)$,  $\partial^2 \mathcal L (., \bar \lambda)(d) $ is locally bounded around $\bar x$. This, together with the fact that  $\xi_k \to \bar x$, implies that $\{z_k\}$ is bounded.  Without loss of generality, we may assume that $z_k$ converges to some $z\in\partial^2 \mathcal{L} (\bar x, \bar \lambda)(d)$. By passing to the limit as $k\to \infty$ in \eqref{New-equa-1}, we obtain $$\dfrac{1}{2} \langle z, d \rangle  \ge 0,$$
	as required.
	
	Now we assume that $\bar x$ is a local optimal solution to the problem \eqref{Problem} and the  locally complementary slackness condition holds at this point with respect to a set $\bar J\in \mathcal{J}(\bar x)$ and a neighborhood $U$ of $\bar x$. Then, by the ACQ,  $\bar x$ is an $L$-stationary point, i.e,
	\begin{equation*}
		0\in\nabla f(\bar x)+\bigcap_{J\in\mathcal{J}(\bar x)}\mathrm{span}\{\nabla g_i(\bar x)\,:\, i\in J^c\}.  
	\end{equation*}
	Hence, for  $\bar J\in \mathcal{J}(\bar x)$ as above, there exists  $\bar \lambda\in S_{\bar J^c}$ such that 
$$\nabla f(\bar x)+\sum_{i\in\bar J^c}\bar\lambda_i\nabla g_i(\bar x)=0,$$  
i.e., $\bar \lambda\in \Lambda(\bar x)$. Now fix any $d\in L(\bar x)$, then by the ACQ,   $d\in T_\Omega(\bar x)$. Hence, we may assume that there exist sequences $t_k\downarrow 0$ and $d_k\to d$ such that $x_k:=\bar x+t_kd_k\in\Omega\cap U$ and $f(x_k)-f(\bar x)\geq 0$ for all $k\in\N$. This, together with the locally complementary slackness condition, implies that $g(x_k)\in S_{\bar J}$ for all $k\in\N$. Thus $\langle \bar\lambda, g(x_k)-g(\bar x)\rangle=0$ for all $k\in\N$. The remainder of the proof proceeds similarly to the previous part.

	\medskip 
	(ii) Since $\bar x$ is an isolated local solution of order $2$ to the problem \eqref{Problem}, there exists $\rho>0$ and a neighborhood $U$ of $\bar x$ such that 
	\begin{equation*}
		f(x)-f(\bar x)\geq \rho \|x-\bar x\|^2 \ \ \forall x\in \Omega\cap U.
	\end{equation*}  
	For $d\in L(\bar x)\setminus\{0\}$, there exist sequences $t_k\downarrow 0$ and $d_k\to d$ such that $x_k:=\bar x+t_kd_k\in\Omega\cap U$ for all $k\in\N$. By a similar argument as in the proof of part (i), we can find $\bar J\in \mathcal{J}(\bar x)$ and $\bar \lambda\in S_{\bar J^c}$ such that
	\begin{equation*}
		\mathcal{L} (x_k, \bar \lambda)- \mathcal{L} (\bar x, \bar \lambda)= f(x_k)-f(\bar x) \ge \rho t_k^2\|d_k\|^2 \ \ \forall k\in\N.
	\end{equation*}
	Once again, by repeating the arguments used in the proof of part (i), we can find $z\in\partial^2 \mathcal{L} (\bar x, \bar \lambda)(d)$  satisfying
	$$\dfrac{1}{2} \langle z, d \rangle  \ge \rho\|d\|^2>0,$$
	as required. 
	
	The proof of the SSSOC property under the the locally complementary slackness condition is similar to that of part (i) and is therefore omitted.
\end{proof}

\begin{remark}\rm The conclusions of  Theorem \ref{Necessary-condition}
	do not remain valid if the statement {\em ``there exists $z\in\partial^2 \mathcal{L} (\bar x, \bar\lambda)(d)$''} is replaced by  {\em``for all $z\in\partial^2\mathcal{L}(\bar x,\bar\lambda)(d)$''}. To illustrate this, consider the following example, which is inspired by \cite[Example 3.2]{ChieuLeeYen2017}. 
\end{remark}

\begin{example}\label{Example-1}\rm  Let $n=3$, $s=2$, $g(x)=x$, $f(x)= \varphi(x_2-1)$ for all $x=(x_1, x_2, x_3)\in\mathbb{R}^3$, where $\varphi\colon\mathbb{R}\to \mathbb{R}$ is the function defined by $\varphi(u):=\int^u_0 h(t)dt$. Here, $h\colon\mathbb{R} \to \mathbb{R}$ is defined as follows:
\medskip
	
	(1) $h(0):=0$;
\medskip
	
	(2) $h(t):=\tfrac{1}{2}t$ for $t>\tfrac{1}{2}$;
\medskip

	(3) if $k\in \mathbb{N}$ and $\tfrac{1}{2^{k+1}}<t\leq\tfrac{3}{2^{k+2}}$, then $h(t):=2t-\tfrac{3}{2^{k+2}}$;
\medskip
	
	(4) if $k\in \mathbb{N}$ and $\tfrac{3}{2^{k+2}}<t\leq \tfrac{1}{2^k}$, then $h(t):=-t+\tfrac{3}{2^{k+1}}$;
\medskip

(5) $h(t)=-h(-t)$ for each $t<0$.

\begin{center}
	\begin{figure}[htp]
		\begin{center}
			\includegraphics[height=5.5cm,width=9cm]{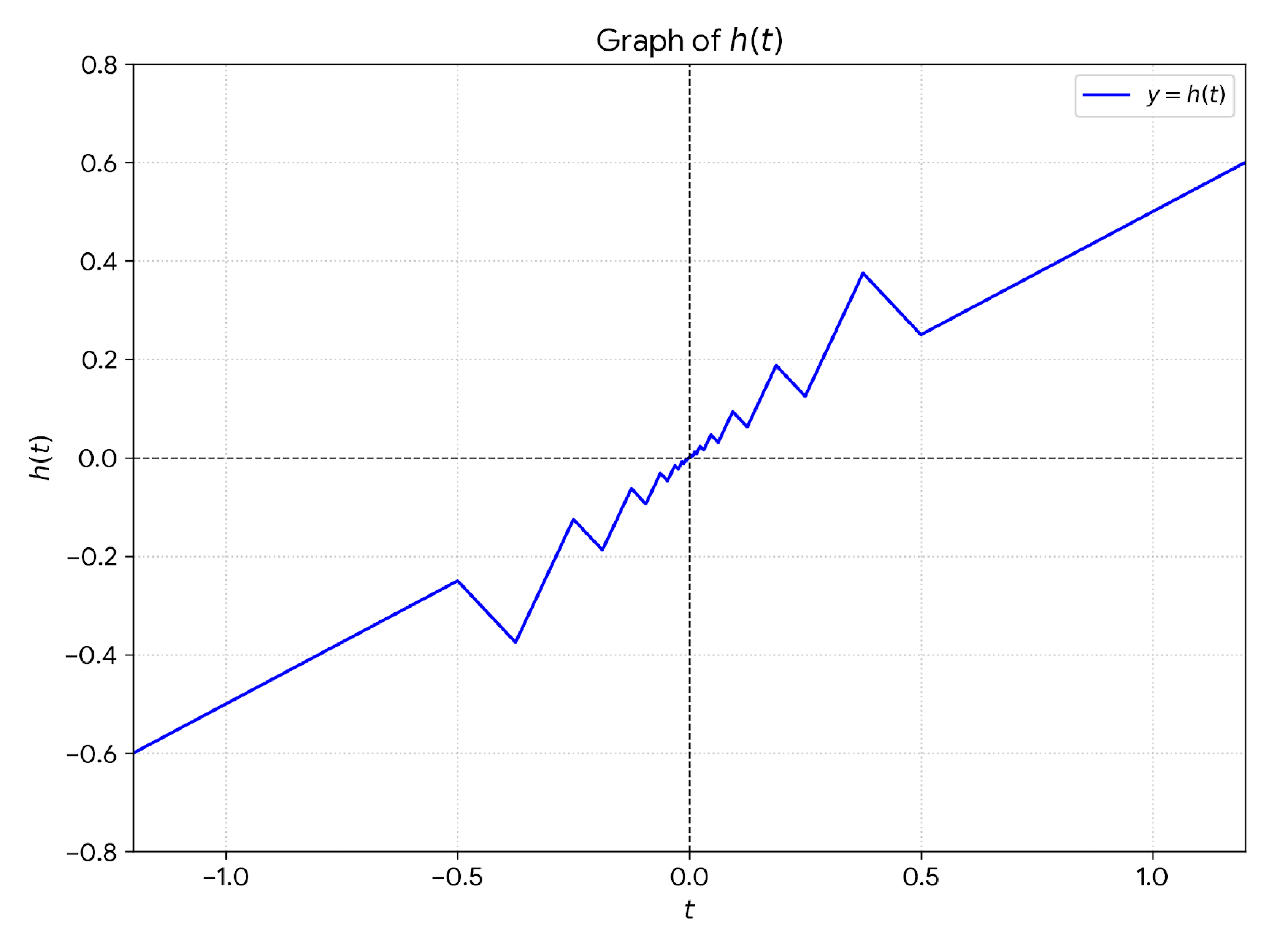}
		\end{center}
		\caption{Graph of  $h$ defined in Example \ref{Example-1}.}\label{Fig_1}
	\end{figure}
\end{center}

Then 
$$\Omega=\mathcal{S}=\{x\in\mathbb{R}^3\,:\, \|x\|_0\leq 2\}.$$
It follows from \cite[Example 3.2]{ChieuLeeYen2017} that $\varphi\in C^{1,1}$, $\bar u=0$ is a global minimizer of $\varphi$, and $1\in \partial^2\varphi(0)(-1)$; see Fig. \ref{Fig_1}.  Thus $\bar x=(1,1, 0)$ is a global minimizer of $f$ on $\Omega$. It is easy to see that   $\Lambda(\bar x)=\{\bar\lambda\}=\{(0,0,0)\}$,  $L(\bar x)=\{d\in \mathbb{R}^3\,:\, d_3=0\}$, and that the MFCQ and the locally complementary slackness condition hold at $\bar x$. Then, by Theorem \ref{Necessary-condition}(i), the SSOC holds at $\bar x$. However, since 
$$\partial^2 \mathcal{L}(\bar x, \bar \lambda)(d)= \{0\}\times\partial^2\varphi(0)(d_2)\times\{0\},$$
we see that $\bar z=(0, 1, 0)\in \partial^2\mathcal{L}(\bar x, \bar \lambda)(\bar d)$, where $\bar d=(0, -1, 0)\in L(\bar x)$, and we therefore get $\langle\bar z, \bar d\rangle <0$.  
\end{example}

\begin{theorem} Let $\bar x\in\Omega$ be a feasible solution to the problem \eqref{Problem} which satisfies the ACQ together with the intersection property \eqref{New-equa-10}.
\begin{enumerate}
	\item [\rm(i)] If $\bar x$ is a local optimal solution to the problem \eqref{Problem}$,$ then $\bar x$ is an $S$-stationary point and the SSOC holds at $\bar x$.
	
	\item [\rm(ii)] If $\bar x$ is an isolated local solution of order $2$ to the problem \eqref{Problem}$,$ then $\bar x$ is an $S$-stationary point and the SSSOC holds at $\bar x$.  
\end{enumerate}	
\end{theorem}
\begin{proof}
By the intersection property and Lemma \ref{Lemma31}(ii), $\bar x$ is an $S$-stationary point. Fix any $\bar \lambda\in \widehat{N}_\mathcal{S}(g(\bar x))$  that satisfies the $S$-stationarity condition:
$$\nabla f(\bar x)+\bar\lambda^T\nabla g(\bar x)=0.$$
Now taking any $d\in L(\bar x)$. Then there exists a  sequence   $x_k\xrightarrow{\Omega} \bar x$ satisfying \eqref{New-equa-11}.  If $\|g(\bar x)\|_0=s$, then, by \eqref{New-equa-11}, $\Gamma(g(x_k))=\Gamma(g(\bar x))$ and $\bar\lambda_i=0$ for all $i\in \Gamma(g(x_k)-g(\bar x))$  and $k\in\mathbb{N}$. This implies that $\langle \bar\lambda, g(x_k)-g(\bar x)\rangle=0$ for all $k\in\N$. If $\|g(\bar x)\|_0<s$, then $\bar\lambda=0$ and so $\langle \bar\lambda, g(x_k)-g(\bar x)\rangle=0$. The remainder of the proof is similar to that of Theorem \ref{Necessary-condition} and is omitted.
\end{proof}

\begin{remark}\rm We note that the intersection property stated in \eqref{New-equa-10} was first introduced in \cite[Definition 9]{Flegel-et-al-07}, whereas the locally complementary slackness condition was introduced in \cite[Theorem 4.2]{Movahedian-et-al}. In fact, these two properties are independent. For example, take $n=m=2$, $s=1$ and $g(x):=(x_1, x_2)$ for all $x=(x_1, x_2)\in\mathbb{R}^2$. Then
$$\Omega=\{x\in\mathbb{R}^2\,:\, \|g(x)\|_0\leq 1\}=\{x\in\mathbb{R}^2\,:\, x_1x_2=0\}.$$	
Let $\bar x=(0,0)\in\Omega$. It is easy to check that the intersection property holds at $\bar x$. However, there do not exist $\bar J\in \mathcal{J}(\bar x)=\{\{1\}, \{2\}\}$ and a neighborhood $U$ of $\bar x$ such that $g(x)\in S_{\bar J}$ for all $x\in\Omega\cap U$, i.e., the locally complementary slackness condition does not hold at $\bar x$.   

Now taking $n=1$, $m=2$, $s=1$, and $g(x):=(x, x)$ for all $x\in\mathbb{R}$, we have
$$\Omega=\{x\in\mathbb{R}\,:\, \|g(x)\|_0\leq 1\}=\{0\}.$$	
Let $\bar x=0\in\Omega$. Clearly, the complementary slackness condition holds at $\bar x$.  However,
$$\bigcap_{J\in \mathcal{J}(\bar x)}\mathrm{span}\{\nabla g_i(\bar x)\,:\, i\in J^c\}=\mathbb{R},$$
i.e., the intersection property does not hold at $\bar x$.    
\end{remark}	
\medskip
 
We now proceed to establish second-order sufficient optimality conditions.
\begin{theorem}[Sufficient optimality conditions]\label{Sufficient-conditions} Assume that $\bar x$ is an $S$-stationary point of \eqref{Problem} with respect to a  Lagrange multiplier vector $\bar \lambda\in\Lambda(\bar x)$. If for any $d\in [T_\Omega(\bar x)\cap\mathrm{ker}\nabla f(\bar x)]\setminus\{0\}$  and $z\in\partial^2 \mathcal{L} (\bar x, \bar \lambda)(d)$ one has $\langle z, d\rangle>0,$ then $\bar x$ is an isolated local solution of order $2$ to the problem \eqref{Problem}.
\end{theorem}
\begin{proof} Suppose to the contrary that $\bar x$ is not an isolated local solution of order $2$ to the problem \eqref{Problem}. Then there exists a sequence $x_k\in \Omega$ such that
	\begin{equation}\label{New-equa-5}
		f(x_k)-f(\bar x)<\frac{1}{k}\|x_k-\bar x\|^2 \ \ \forall k\in\N.
	\end{equation}
	Clearly, $x_k\neq \bar x$ for all $k\in\N$. Hence $t_k:=\|x_k-\bar x\|\downarrow 0$ as $k\to\infty$. Put $d_k:=\frac{x_k-\bar x}{t_k}$, then $\|d_k\|=1$ and so we may assume that $d_k$ converges to some $d\in\R^n$ with $\|d\|=1$. Clearly, $d\in T_\Omega(\bar x)$. By a similar argument as in the proof of Theorem \ref{Necessary-condition}(i), we can find a set $\bar J\in \mathcal{J}(\bar x)$ such that
	\begin{equation}\label{New-equa-3}
		\Gamma(g(\bar x))\cup \Gamma(g(x_k)-g(\bar x))\subset \bar J\in\mathcal{J}(\bar x) \  \ \forall k\in\N.
	\end{equation} 
We now consider two cases of $\bar x$ as follows.
	
	{\bf Case 1.} $\Gamma(g(\bar x))=s$. Then it follows from \eqref{New-equa-3} that $\Gamma(g(x_k)-g(\bar x))=0$, i.e., $g(x_k)=g(\bar x)$ for all $k\in\N$. Hence
	\begin{equation}\label{New-equa-4}
		\mathcal{L}(x_k,\bar\lambda)-\mathcal{L}(\bar x,\bar\lambda)=f(x_k)-f(\bar x) \ \ \forall k\in\N.
	\end{equation}
	
	{\bf Case 2.} $\Gamma(g(\bar x))<s$. Since $\bar x$ is an $S$-stationary point, $\bar\lambda=0$ and we also get \eqref{New-equa-4}.
	
	We now show that $\langle \nabla f(\bar x), d\rangle=0$. Indeed, by \eqref{New-equa-4} and the mean value theorem, we can find $u_k$ and $\xi_k$  belonging to $[\bar x, x_k]$ such that
	\begin{equation*}
		\langle\nabla_x\mathcal{L}(u_k, \bar \lambda), d_k\rangle=\langle\nabla f(\xi_k), d_k\rangle \ \ \forall k\in\N.
	\end{equation*}
	Hence, by letting $k\to\infty$, we obtain
	\begin{equation*}
		0= \langle\nabla_x\mathcal{L}(\bar x, \bar \lambda), d\rangle=\langle\nabla f(\bar x), d\rangle,
	\end{equation*}
	as required.
	
	It follows from \eqref{New-equa-5} that
	\begin{equation*}
		\mathcal{L}(x_k,\bar\lambda)-\mathcal{L}(\bar x,\bar\lambda)\leq \frac{1}{k}t_k^2\|d_k\|^2 \ \ \forall k\in\N.
	\end{equation*}
	
	Following the proof of Theorem \ref{Necessary-condition}(i), we write $\mathcal{L}(x_k,\bar\lambda)-\mathcal{L}(\bar x,\bar\lambda)$ in the form $$\mathcal{L}(x_k,\bar\lambda)-\mathcal{L}(\bar x,\bar\lambda)=\mathcal{L}_{1k}+\mathcal{L}_{2k},$$   
	where $\mathcal{L}_{1k}:=	\mathcal{L}(x_k,\bar\lambda)  -\mathcal{L}(\bar x+t_k d, \bar\lambda)$ and 
	$$\mathcal{L}_{2k}:=\mathcal{L} (\bar x + t_k d, \bar \lambda)- \mathcal{L} (\bar x, \bar \lambda) -  \langle \nabla_x \mathcal{L}(\bar x, \bar \lambda) , t_k d \rangle.$$
	Then $\frac{\mathcal{L}_{1k}}{t_k^2}\to 0$ as $k\to\infty$ due to \eqref{New-equa-7}. Furthermore,  by applying Taylor formula (Theorem \ref{Taylor_formula}) and Proposition \ref{property}  to the function $\mathcal{L}(\cdot, \bar\lambda)$ on the interval $[\bar x, \bar x + t_k d]$, we obtain a sequence  $z_k\in \partial^2\mathcal{L}(\xi_k,\bar \lambda)(d)$, where $\xi_k\in [\bar x, \bar x + t_k d]$, such that $z_k\to z$ for some $z\in \partial^2\mathcal{L}(\bar x,\bar \lambda)(d)$ and 
	\begin{equation*}
		\frac{1}{2}t_k^2\langle z_k, d\rangle\leq \mathcal{L}_{2k} \  \ \forall k\in\N.  
	\end{equation*}
	Hence,
	\begin{equation}\label{New-equa-6}
		\mathcal{L}_{1k}+\frac{1}{2}t_k^2\langle z_k, d\rangle\leq \mathcal{L}(x_k,\bar\lambda)-\mathcal{L}(\bar x,\bar\lambda)\leq \frac{1}{k}t_k^2\|d_k\|^2 \ \ \forall k\in\N.
	\end{equation}
	After dividing both sides of \eqref{New-equa-6} by $t_k^2$ and passing to the limit as $k\to\infty$, we can find $z\in \partial^2\mathcal{L}(\bar x,\bar \lambda)(d)$ satisfying  $\langle z, d\rangle\leq 0$,  which contradicts the assumptions.
\end{proof}

\section{Sparse multiobjective optimization problems}\label{Section 4}
In this section, we consider the following sparse multiobjective optimization problem
\begin{equation}\label{MOP}
	\mathrm{Min}\,_{\mathbb{R}^q_+}F(x):=(f_1(x), \ldots, f_q(x))\ \ \text{s.t.}\ \ x\in\Omega,\tag{SMP}
\end{equation}
where $\Omega$ is defined as in \eqref{Problem} and  $f_l\in C^{1,1}(\mathbb{R}^n)$, $l\in L:=\{1, \ldots, q\}$. 

We now define the Lagrangian function associated with the problem~\eqref{MOP} by
\begin{align*}
	\mathcal{L}^M (x, \alpha, \lambda)=\alpha^TF(x)+\lambda^Tg(x) \ \ \forall (x, \alpha, \lambda)\in \R^n\times\R^q\times\R^m.
\end{align*}

\begin{definition}\rm  Let $\bar x\in \Omega$. The set of {\em Lagrange multipliers} of \eqref{MOP} at $\bar x$ is defined by
	\begin{equation*}
		\Lambda^M(\bar x):=\{(\alpha,\lambda)\in (\mathbb{R}^q_+\setminus\{0\})\times\mathbb{R}^m\,:\, \exists J\in \mathcal{J}(\bar x) \ \ \text{with}\ \  \lambda \in S_{J^c}, \nabla_x\mathcal{L}^M(\bar x, \alpha, \lambda)=0\}.
	\end{equation*}	
\end{definition}

For $\bar x\in\Omega$, we denote the {\em strict linearized cone} to $\Omega$ at $\bar x$ by
\begin{equation*}
	L^s(\bar x)=\bigcap_{J\in\mathcal{J}(\bar x)}\{d\in\R^n\,:\, \langle\nabla g_i(\bar x), d\rangle=0 \ \ \forall i\in J^c\}.
\end{equation*}

\begin{definition}[{see, e.g. \cite{Ginchev_etal}}]\rm Let $\bar x\in \Omega$. We say that:
	\begin{enumerate}
		\item[\rm(i)] $\bar x$ is a {\em local weakly Pareto efficient solution} of \eqref{MOP} if there exists a neighborhood $U$ of $\bar x$ such that  
		\begin{equation*}
			\max\{f_1(x)-f_1(\bar x), \ldots, f_q(x)-f_q(\bar x)\}\geq 0 \ \ \forall x\in \Omega\cap U. 
		\end{equation*}
		
		\item[\rm(ii)] 	$\bar x$ is a {\em local strictly Pareto efficient solution} of \eqref{MOP} if there exists a neighborhood $U$ of $\bar x$ such that  
		\begin{equation*}
			\min\{f_1(x)-f_1(\bar x), \ldots, f_q(x)-f_q(\bar x)\}\geq 0 \ \ \forall x\in \Omega\cap U. 
		\end{equation*}
		
		\item[\rm(iii)]  $\bar x$  is an {\em  isolated local solution of order $2$} of \eqref{MOP} if there exist a positive number $\beta$ and a neighborhood $U$ of $\bar x$ such that
		\begin{equation*}
			\max\{f_1(x)-f_1(\bar x), \ldots, f_q(x)-f_q(\bar x)\}\geq \beta \|x-\bar x\|^2 \ \ \forall x\in \Omega\cap U.
		\end{equation*}
	\end{enumerate}
\end{definition}

\begin{definition}\rm Let $\bar x\in\Omega$. We say that
	\begin{enumerate}
		\item[\rm(i)] $\bar x$ is an {\em  $M$-stationary point} of \eqref{MOP} if
		\begin{equation*}
			0\in\mathrm{co}\,\{\nabla f_l(\bar x)\,:\, l\in L\}+\bigcup_{J\in \mathcal{J}(\bar x)}\mathrm{span}\{\nabla g_i(\bar x)\,:\, i\in J^c\}.
		\end{equation*}
		\item[\rm(ii)] $\bar x$ is an {\em  $L$-stationary point} of \eqref{MOP} if
		\begin{equation*}
			0\in \mathrm{co}\,\{\nabla f_l(\bar x)\,:\, l\in L\}+\bigcap_{J\in \mathcal{J}(\bar x)}\mathrm{span}\{\nabla g_i(\bar x)\,:\, i\in J^c\}.
		\end{equation*}
		\item[\rm(iii)] $\bar x$ is an {\em  $S$-stationary point} of \eqref{MOP} if
		\begin{equation*}
			\begin{cases}
				0\in\mathrm{co}\,\{\nabla f_l(\bar x)\,:\, l\in L\}+\mathrm{span}\{\nabla g_i(\bar x)\,:\, i\notin \Gamma(g(\bar x))\},  \ \ &\text{if}\ \ |\Gamma(g(\bar x))|=s,
				\\
				0\in\mathrm{co}\,\{\nabla f_l(\bar x)\,:\, l\in L\}, \ \ &\text{otherwise}.
			\end{cases}
		\end{equation*}
	\end{enumerate}
\end{definition}

The following theorem provides a first-order necessary optimality condition  for a local strictly Pareto efficient solution to the problem \eqref{MOP}.   
\begin{proposition}\label{1-Necessary-condition-M} Let $\bar x\in\Omega$ be a feasible solution to the problem \eqref{MOP} which satisfies the GCQ. If $\bar x$ is a local strictly Pareto efficient solution to the problem \eqref{MOP}$,$ then $\bar x$ is an $L$-stationary point of \eqref{MOP}. 	In other words$,$ $\Lambda^M(\bar x)$ is a nonempty set.
\end{proposition}
\begin{proof} Since $\bar x$ is a local strictly Pareto efficient solution to the problem \eqref{MOP}, $\bar x$ is a local minimizer of the function $f_l(\cdot)-f_l(\bar x)$ on $\Omega$ for all $l\in L$. By Lemma \ref{Lemma31}, one has
	\begin{equation*}
		0\in\nabla f_l(\bar x)+\bigcap_{J\in\mathcal{J}(\bar x)}\mathrm{span}\{\nabla g_i(\bar x)\,:\, i\in J^c\} \ \ \forall l\in L.
	\end{equation*}     
	Hence,	 
	\begin{equation*}
		0\in  \mathrm{co}\,\{\nabla f_l(\bar x)\,:\, l\in L\} +\bigcap_{J\in\mathcal{J}(\bar x)}\mathrm{span}\{\nabla g_i(\bar x)\,:\, i\in J^c\}.
	\end{equation*}
	This means that $\Lambda^M(\bar x)$ is nonempty. 
\end{proof}

We now establish second-order necessary optimality conditions of \eqref{MOP}. 
\begin{theorem} \label{Necessary-condition-M} Let $\bar x\in\Omega$ be a feasible solution to the problem \eqref{MOP} which satisfies the ACQ. 
	\begin{enumerate}
		\item[\rm(i)] If $\bar x$ is a local weakly Pareto efficient solution to the problem \eqref{MOP}$,$ then for each $d\in L^s(\bar x)$, there exist $\bar J\in\mathcal{J}(\bar x)$ and $(\bar \alpha, \bar\lambda)\in (\R^q_+\setminus\{0\})\times S_{\bar J^c}$ such that either
		\begin{equation*}
			\langle\nabla_x\mathcal{L}^M(\bar x, \bar\alpha,\bar \lambda), d\rangle>0,
		\end{equation*} 
		or
		\begin{equation*}
			\langle\nabla_x\mathcal{L}^M(\bar x, \bar\alpha, \bar\lambda), d\rangle=0 \ \ \text{and}\ \ \langle z, d\rangle\geq 0 \ \ \text{for some}\ \ z\in \partial^2\mathcal{L}^M(\bar x,\bar \alpha, \bar\lambda)(d).  
		\end{equation*}
		
		\item[\rm(ii)] If $\bar x$ is an isolated local solution of order $2$ to the problem \eqref{MOP}$,$ then for each $d\in L^s(\bar x),$ there exist $\bar J\in\mathcal{J}(\bar x)$ and $(\bar \alpha, \bar\lambda)\in (\R^q_+\setminus\{0\})\times S_{\bar J^c}$ such that either
		\begin{equation*}
			\langle\nabla_x\mathcal{L}^M(\bar x, \bar\alpha,\bar \lambda), d\rangle>0,
		\end{equation*} 
		or
		\begin{equation*}
			\langle\nabla_x\mathcal{L}^M(\bar x, \bar\alpha, \bar\lambda), d\rangle=0 \ \ \text{and}\ \ \langle z, d\rangle> 0 \ \ \text{for some}\ \ z\in \partial^2\mathcal{L}^M(\bar x,\bar \alpha, \bar\lambda)(d).  
		\end{equation*}
		\item[\rm(iii)] If $\bar x$ is a local strictly Pareto efficient solution to the problem \eqref{MOP}$,$ then for each $d\in L(\bar x),$ there exist  $(\bar \alpha, \bar \lambda)\in\Lambda^M(\bar x)$ and $z\in \partial^2\mathcal{L}^M(\bar x,\bar \alpha, \bar\lambda)(d)$ such that $\langle z, d\rangle\geq 0$.
	\end{enumerate}	 
\end{theorem}
\begin{proof} (i) Let $d\in L^s(\bar x)$, then $d\in T_\Omega(\bar x)$. Hence there exist sequences $t_k\downarrow 0$ and $d_k\to d$ as $k\to\infty$ such that $x_k:=\bar x+t_kd_k\in\Omega$ for all $k\in\N$. Since $\bar x$ is a local weakly Pareto efficient solution to the problem \eqref{MOP}, there exists a neighborhood $U$ of $\bar x$ such that
	\begin{equation*}
		\max\{f_l(x)-f_l(\bar x)\,:\, l\in L\}\geq 0 \ \ \forall x\in\Omega\cap U.
	\end{equation*}
	Hence, we may assume that
	\begin{equation*}
		\max\{f_l(x_k)-f_l(\bar x)\,:\, l\in L\}\geq 0 \ \ \forall k\in\N.
	\end{equation*}
	This implies that there exists $\bar \alpha\in\R^q_+\setminus\{0\}$ such that
	\begin{equation}\label{New-equa-8}
		\sum_{l=1}^{q}\bar\alpha_l[f_l(x_k)-f_l(\bar x)]\geq 0
	\end{equation}	
	for infinitely many $k\in\N$. By passing to a subsequence if necessary, we may assume that \eqref{New-equa-8}  holds for all $k\in\N$. Hence, for each $k\in\N$ one has
	\begin{equation*}
		\sum_{l=1}^{q}\bar\alpha_l\langle\nabla f_l(\xi_{lk}), d_k\rangle\geq 0 \ \ \text{for some}\ \ \xi_{lk}\in [\bar x, x_k].
	\end{equation*}  
	Letting $k\to\infty$ we get  $\sum_{l=1}^{q}\bar\alpha_l\langle\nabla f_l(\bar x), d\rangle\geq 0$. Now, by a similar argument as in the proof of Theorem \ref{Necessary-condition}(i), we can find a set $\bar J\in \mathcal{J}(\bar x)$ and  such that
	\begin{equation*} 
		\Gamma(g(\bar x))\cup \Gamma(g(x_k)-g(\bar x))\subset \Gamma(g(x_k)) \subset \bar J\in\mathcal{J}(\bar x) \  \ \forall k\in\N.	
	\end{equation*} 
	Let $\bar \lambda\in S_{\bar J}^c$, then $\langle\bar\lambda, g(x_k)-g(\bar x)\rangle=0$ for all $k\in\N$. On the other hand, since $d\in L^s(\bar x)$, one has $\langle\nabla g_i(\bar x), d\rangle=0$ for all $i\in \bar J^c$. Hence
	\begin{equation*}
		\langle\nabla_x\mathcal{L}^M(\bar x, \bar\alpha,\bar \lambda), d\rangle= \sum_{l=1}^{q}\bar\alpha_l\langle\nabla f_l(\bar x), d\rangle+\sum_{i\in\bar J^c}\bar \lambda_i\langle\nabla g_i(\bar x), d\rangle=\sum_{l=1}^{q}\bar\alpha_l\langle\nabla f_l(\bar x), d\rangle\geq 0
	\end{equation*}    
	and
	\begin{equation*}
		\mathcal{L}^M(x_k,\bar\alpha, \bar\lambda)-\mathcal{L}^M(\bar x,\bar\alpha, \bar\lambda)=\sum_{l=1}^{q}\bar\alpha_l[f_l(x_k)-f_l(\bar x)]\geq 0.
	\end{equation*}
	If $\langle\nabla_x\mathcal{L}^M(\bar x, \bar\alpha,\bar \lambda), d\rangle=0$, then an argument similar to that used in the proof of Theorem \ref{Necessary-condition}(i) shows that there exists  $z\in \partial^2\mathcal{L}^M(\bar x,\bar \alpha, \bar\lambda)(d)$  such that $\langle z, d\rangle\geq 0$, as required.
	
	(ii) The proof is quite similar to that of part (i) and is therefore omitted.
	
	(iii) For $d\in L(\bar x)$, let $t_k\downarrow 0$ and $d_k\to d$ be such that $x_k:=\bar x+t_kd_k\in\Omega$ for all $k\in\N$. Since $\bar x$ is a local strictly Pareto efficient solution to the problem \eqref{MOP}, there exists a neighborhood $U$ of $\bar x$ such that
	\begin{equation*}
		f_l(x)-f_l(\bar x)\geq 0 \ \ \forall x\in\Omega\cap U, \ \  \forall l\in L.
	\end{equation*} 
	Let $\bar J\in\mathcal{J}(\bar x)$ and $\bar \lambda\in S_{\bar J}^c$ be such that $\langle\bar\lambda, g(x_k)-g(\bar x)\rangle=0$ for all $k\in\N$. Since $\bar x$ is a local strictly Pareto efficient solution to the problem \eqref{MOP},  $\Lambda^M(\bar x)$  is nonempty due to  Proposition \ref{1-Necessary-condition-M}. Let $\bar\alpha\in\R^q_+\setminus\{0\}$ be such that $\nabla_x\mathcal{L}^M(\bar x, \bar\alpha,\bar \lambda)=0$. Then 
	\begin{equation*}
		\mathcal{L}^M(x_k,\bar\alpha, \bar\lambda)-\mathcal{L}^M(\bar x,\bar\alpha, \bar\lambda)=\sum_{l=1}^{q}\bar\alpha_l[f_l(x_k)-f_l(\bar x)]\geq 0
	\end{equation*} 
	for all $k$ large enough. Analysis similar to that in the proof of part (i), there exists $z\in \partial^2\mathcal{L}^M(\bar x,\bar \alpha, \bar\lambda)(d)$ such that $\langle z, d\rangle\geq 0$. The proof is complete. 
\end{proof}
 
To establish sufficient optimality conditions for \eqref{MOP}, we introduce the following critical cone.
\begin{definition}\rm  Let $\bar x\in\Omega$. The {\em cone of  critical directions} of \eqref{MOP} at $\bar x$ is defined by
	\begin{equation*}
		\mathcal{K}(\bar x):=\{d\in T_\Omega(\bar x)\,:\, \langle\nabla f_l(\bar x), d\rangle\leq 0 \ \forall l\in L, \langle\nabla f_k(\bar x), d\rangle=0  \ \ \text{for some}\ \ k\in L\}.
	\end{equation*}
	
\end{definition}

We conclude this section by presenting a second-order sufficient optimality condition for an isolated local solution of order $2$ to the problem \eqref{MOP}.

\begin{theorem} Assume that $\bar x\in \Omega$ is a S-stationary point and $(\bar \alpha, \bar \lambda)\in\Lambda^M(\bar x)$. If for every $d\in \mathcal{K}(\bar x)\setminus\{0\}$ and $z\in\partial^2\mathcal{L}^M(\bar x, \bar\alpha, \bar \lambda)$ one has
	\begin{equation*}
		\langle z, d\rangle>0,
	\end{equation*}
	then $\bar x$ is  an isolated local solution of order $2$ to the problem \eqref{MOP}.
\end{theorem}
	\begin{proof} Suppose to the contrary that $\bar x$ is not an isolated local solution of order $2$ of~\eqref{MOP}. Then, by definition, there exists a sequence $x_k \xrightarrow{\Omega\setminus\{\bar x\}} \bar x$ such that  
		\begin{align}\label{xk-MOP}
			\max\{f_l(x_k)-f_l(\bar x)\,:\, l\in L\}\leq \frac{1}{k}\|x_k-\bar x\|^2, \ \forall k\in \mathbb{N}.
		\end{align}
		Put	$t_k:=\| x_k-\bar x\|$ and $d_k:=\frac{x_k-\bar x}{t_k}$. Then, $t_k\downarrow 0$ and we may assume that $d_k\to d$ as $k\to\infty$. Clearly, $d\in T_\Omega(\bar x)$ and $\|d\|=1$. By a similar argument as in the proof of Theorem \ref{Sufficient-conditions}, we see that
		\begin{equation}\label{New-equa-9}
			\mathcal{L}^M(x_k, \bar\alpha, \bar\lambda)-\mathcal{L}^M(\bar x, \bar\alpha, \bar\lambda)=\sum_{l=1}^q\bar\alpha_l (f_l(x_k)-f_l(\bar x))  \ \ \forall k\in\N.
		\end{equation}
		It follows from \eqref{xk-MOP} that
		\begin{equation*}
			f_l(x_k)-f_l(\bar x)\leq \frac{1}{k}t_k^2\|d_k\|^2, \ \ \forall l\in L, k\in \mathbb{N}.
		\end{equation*} 
		By the  mean value theorem, we can find $\xi_k\in[\bar x, x_k]$ such that
		\begin{equation*}
			\langle\nabla f_l(\xi_k), d_k\rangle=\frac{f_l(x_k)-f_l(\bar x)}{t_k}	 \leq \frac{1}{k}t_k\|d_k\|^2.  
		\end{equation*} 
		Hence, by letting $k\to\infty$, we obtain $\langle\nabla f_l(\bar x), d\rangle\leq 0$ for all $l\in L$. We claim that there exists $l\in L$ such that $\langle\nabla f_l(\bar x), d\rangle= 0$, or, equivalently, $d\in\mathcal{K}(\bar x)$. If this is not case,  then 
		$$\sum_{l=1}^q\bar\alpha_l\langle \nabla f_l(\bar x), d\rangle<0.$$
		By \eqref{New-equa-9} and the mean value theorem we can find $u_k$ and $\xi_{lk}$ belonging to $[\bar x, x_k]$ such that 
		\begin{equation*}
			\nabla_x\mathcal{L}^M(u_k,\bar\alpha, \bar\lambda)=\sum_{l=1}^q\bar\alpha_l\langle \nabla f_l(\xi_{lk}), d\rangle.
		\end{equation*} 
		Letting $k\to\infty$ we obtain
		\begin{equation*}
			0=\nabla_x\mathcal{L}^M(\bar x,\bar\alpha, \bar\lambda)=\sum_{l=1}^q\bar\alpha_l\langle \nabla f_l(\bar x, d\rangle <0,
		\end{equation*}
		a contradiction. Hence $d\in\mathcal{K}(\bar x)\setminus\{0\}$. 
		
		We now write $\mathcal{L}^M(x_k, \bar\alpha, \bar\lambda)-\mathcal{L}^M(\bar x, \bar\alpha, \bar\lambda)$ as
		\begin{equation*}
			\mathcal{L}^M(x_k, \bar\alpha, \bar\lambda)-\mathcal{L}^M(\bar x, \bar\alpha, \bar\lambda)= \mathcal{L}^M_{1k}+\mathcal{L}^M_{2k},
		\end{equation*}
		where $\mathcal{L}^M_{1k}= \mathcal{L}^M(x_k, \bar\alpha, \bar\lambda)- \mathcal{L}^M(\bar x+t_kd, \bar\alpha, \bar\lambda)$ and 
		\begin{equation*}
			\mathcal{L}^M_{2k}= \mathcal{L}^M(\bar x+t_kd, \bar\alpha, \bar\lambda) - \mathcal{L}^M(\bar x, \bar\alpha, \bar\lambda)-\langle\nabla_x\mathcal{L}^M(\bar x,\bar\alpha, \bar\lambda), t_kd \rangle. 
		\end{equation*}
		Analysis similar to that in the proof of Theorem \ref{Sufficient-conditions}, we obtain that $\frac{\mathcal{L}^M_{1k}}{t_k^2}\to 0$ and
		\begin{equation*}
			\frac{\mathcal{L}^M_{1k}}{t_k^2}+\frac{1}{2}\langle z_k, d\rangle\leq \mathcal{L}^M(x_k, \bar\alpha, \bar\lambda)-\mathcal{L}^M(\bar x, \bar\alpha, \bar\lambda)\leq \frac{1}{k}\|d_k\|^2\bigg(\sum_{l=1}^q\bar\alpha_l\bigg),
		\end{equation*}
		for some $z_k\in \partial^2\mathcal{L}^M(u_k, \bar\alpha, \bar \lambda)(d)$ and $u_k\in [\bar x, \bar x+t_kd]$ for all $k\in\N$. By letting $k\to\infty$, we get $\langle z, d\rangle\leq 0$ for some $z\in \partial^2\mathcal{L}^M(\bar x, \bar\alpha, \bar \lambda)(d)$, a contradiction. The proof is complete. 
	\end{proof}


\section{Conclusions}\label{Conclusions}
 
In this paper, we have investigated second-order optimality conditions for a class of $C^{1,1}$ sparsity constrained optimization problems by exploiting their disjunctive structure. The sparsity constraint was represented as a finite union of linear subspaces, which enabled us to employ tools from variational analysis and disjunctive programming to derive second-order optimality conditions under relatively mild assumptions.

A central feature of our approach is the use of the limiting   second-order subdifferential of the associated Lagrangian function. In contrast to the Fr\'echet second-order subdifferential, which may be empty even for $C^{1,1}$ functions, the limiting second-order subdifferential provides a robust framework for second-order analysis in the present setting. Based on this framework, we established second-order necessary optimality conditions for local optimal solutions and derived second-order sufficient conditions under appropriate assumptions. The obtained results complement and extend several existing second-order optimality conditions for sparsity constrained optimization problems.

We have also applied the developed theory to sparse multiobjective optimization problems. In particular, second-order necessary and sufficient optimality conditions for efficient solutions were obtained by combining the disjunctive structure of the sparsity constraint with the limiting second-order analysis developed in this paper.  

There are several directions for future research. One possible direction is to investigate second-order optimality conditions for more general sparsity and cardinality-constrained problems, including problems with additional nonlinear constraints and more general disjunctive constraint sets. Another interesting topic is to develop corresponding stability and sensitivity results for sparse optimization problems. It would also be worthwhile to study second-order conditions under weaker constraint qualifications and to investigate algorithmic consequences of the proposed conditions. Finally, extending the present approach to broader classes of nonsmooth or merely $C^1$ functions, as well as to other classes of sparse multiobjective optimization problems, remains an interesting subject for further study.

\subsection*{Acknowledgments}  
A part of this work was done while the third author was visiting Academy for Advanced Interdisciplinary Studies, Northeast Normal University, China in April--May, 2026. He is grateful to the Academy for its hospitality and support.

\subsection*{Funding} The research of L.T.T. Huyen and N.V. Tuyen is funded by Hanoi Pedagogical University 2 Foundation for Sciences and Technology Development under Grant number HPU2.2025-UT-10. The research of L.G. Jiao is partially supported by the Chinese National Natural Science Foundation under grant number 12471478.

\section*{Declarations} 
{\bf Data availability} The authors confirm that the data supporting the findings of this study are available within the article.

\medskip
\noindent{\bf Conflict of interest} The authors declare that they have no conflict of interest.

\bibliographystyle{amsplain}

\end{document}